\documentclass{article}
\usepackage{amsmath}
\usepackage{amssymb}
\usepackage{amscd}

\def\no{\noindent}
\numberwithin{equation}{section}

\newtheorem{lemma}[equation]{Lemma}

\newtheorem{theorem}[equation]{Theorem}

\newcommand{\ra}{\rightarrow}
\newcommand{\restr}{\mbox{\Large \(|\)\normalsize}}
\newcommand{\N}{\mathbb N}
\newcommand{\R}{\mathbb R}

\newcommand{\E}{\mathbb E}
\newcommand{\Z}{\mathbb Z}

\newcommand{\qed}{{\unskip\nobreak\hfil
        \penalty50\hskip1em\hbox{}\nobreak\hfil
        $\square$\parfillskip=0pt\finalhyphendemerits=0 \par}}
\newcommand{\proof}{\no{\em Proof.\ }}

\newcommand{\al}{\alpha}

\def\eps{\epsilon}

\def\ol{\overline}

\def\ra{\rightarrow}

\def\defeq{:=}
\begin{document}

\title{Separated nets in Euclidean space and Jacobians of
biLipschitz maps}
\author{Dmitri Burago \thanks{Supported by a Sloan Foundation
Fellowship and NSF grant
DMS-95-05175.}\\
Bruce Kleiner\thanks {Supported by a Sloan Foundation
Fellowship and NSF grants
DMS-95-05175 and DMS-96-26911.}}
\date{\today}
\maketitle

\begin{abstract}
We show that there are separated nets in the Euclidean plane 
which are not biLipschitz  equivalent to the
integer 
lattice. The argument is based on the construction of a continuous
function which is not the  Jacobian of a biLipschitz map. 
\end{abstract}

\section{Introduction}
A subset $X$ of a metric space $Z$ is a {\em separated net}
if there are constants $a,\,b>0$ such that $d(x,x')>a$ for every 
pair $x,\,x'\in X$, and $d(z,X)<b$ for every $z\in Z$.  Every
metric space contains separated nets: they may be constructed
by finding maximal subsets with the property that all pairs of 
points are separated by some distance $a>0$.
It follows easily from the definitions that two spaces are quasi-isometric if 
and only if they contain biLipschitz equivalent separated nets.
One may ask if the choice of these nets matters, or, in other
words,
whether any two separated nets in a given space are biLipschitz
equivalent.  To the best of our knowledge, this problem was first posed by
Gromov \cite[p.23]{asyinv}. 
  The answer is known to be yes for separated nets in non-amenable spaces
 (under mild assumptions about local geometry), see \cite{Gromov,McMullen,Whyte};
more constructive proofs in the case of trees or hyperbolic
groups can be found in \cite{Papas,Bogop}.

In this paper, we prove the following theorem:

\begin{theorem}
\label{T1}
There exists a separated net in the Euclidean plane 
which is not biLipschitz equivalent to the integer lattice.
\end{theorem}

The proof of Theorem \ref{T1} is based on the following result:

\begin{theorem}
\label{T2}
Let $I:=[0,1]$.  Given $c>0$, there is a continuous function
 $\rho:I^2\ra [1,1+c]$,
such that there is no biLipschitz 
map $f:I^2\ra  \E^2$ with 

$$Jac(f):=Det(Df)=\rho \quad a.e.$$

\end{theorem}

\no
{\bf Remarks}

1.  Although we formulate and prove these theorems in the
2-dimensional case, 
the same  proofs work with minor modifications in
higher dimensional  Euclidean spaces as well.
We  only consider the 2-dimensional case here to avoid cumbersome notation.

2. Theorem \ref{T2} also works for Lipschitz homeomorphisms; we
do not use the lower Lipschitz bound on $f$.

3. After the first version of this paper had been written,
Curt McMullen informed us that he also had a proof of Theorems
\ref{T1} and \ref{T2}.  See \cite{McMullen} for a discussion
of the the linear analog of Theorem \ref{T2},
and the H\"older analogs of the mapping problems
in Theorems \ref{T1} and \ref{T2}.

The problem of prescribing Jacobians of homeomorpisms has been
studied by several authors.  In \cite{DacMos} 
Dacorogna and Moser proved that every
$\al$-Holder continuous function is locally the Jacobian of 
a $C^{1,\al}$ homeomorphism, and they then raised the question 
of whether any continuous function is (locally) the
Jacobian of a $C^1$ diffeomorphism.  \cite{RivYe,Ye} consider
the prescribed Jacobian problem in other regularity classes, including the
cases when the Jacobian is in $L^{\infty}$ or in a Sobolev
space.

\bigskip
\centerline{\large Overview of the proofs}

{\em Theorem \ref{T2} implies Theorem \ref{T1}.} 
 Let $\rho:I^2\ra\R$ be measurable with $0<\inf\rho\leq \sup\rho <\infty$.
 We will indicate why $\rho$
 would be the Jacobian of a biLipschitz map $f:I^2\ra\E^2$ if
 all separated nets in $\E^2$ were
biLipschitz equivalent.  
Take a disjoint collection of
squares $S_i\subset\E^2$ with side lengths $l_i$ tending to infinity,
and
``transplant'' $\rho$ to each $S_i$ using appropriate
similarities $\al_i:I^2\ra S_i$, i.e. set $\rho_i\defeq \rho\circ\al_i^{-1}$.  Then
construct
a separated net $L\subset \E^2$ so that the ``local average density'' 
of $L$ in each square $S_i$ approximates $\rho_i^{-1}$.
If $g:L\ra \Z^2$ is a biLipschitz homeomorphism, 
consider ``pullbacks'' of $g\restr_{S_i}$  
to $I^2$, i.e.  pre and post-compose $g\restr_{S_i}$
with suitable similarities so as to get a sequence
of uniformly biLipschitz maps $g_i:I^2\supset Z_i\ra\E^2$.
Then extract a convergent subsequence of the $g_i$'s
via the Arzela-Ascoli theorem, and obtain a limit map 
$f:I^2\ra\E^2$ with  
Jacobian $\rho$.

{\em Theorem \ref{T2}.}
The observation underlying our construction is that
if the Jacobian of $f:I^2\ra\E^2$ oscillates in a 
rectangular neighborhood $U$ of a segment $\ol{xy}\subset I^2$,
then $f$ will be forced to stretch for one of two
reasons: either it maps $\ol{xy}$ to a curve which 
is far from a geodesic between its endpoints, or it
maps $\ol{xy}$ close to the segment $\ol{f(x)f(y)}$
but it sends $U$ to a neighborhood of
$\ol{f(x)f(y)}$ with wiggly boundary in order to have the
correct Jacobian.  By arranging that $Jac(f)$ 
oscillates in neighborhoods of a hierarchy of
smaller and smaller segments we can force $f$ to stretch 
more and more at smaller and smaller scales, eventually
contradicting the Lipschitz condition on $f$.

We now give a more detailed sketch of the proof.  
 
  We first observe that
it is enough to construct, for every $L>1,\,\bar c>0$, a 
continuous function $\rho_{L,\bar c}:I^2\ra[1,1+\bar c]$
such that $\rho_{L,\bar c}$ is not the Jacobian of an 
L-biLipschitz map $I^2\ra\E^2$.  Given such a family
of functions, we can build a new continuous function 
$\rho:I^2\ra[1,1+c]$ which is not the Jacobian 
of any biLipschitz map $I^2\ra\E^2$ as follows.
Take a sequence of disjoint squares $S_k\subset I^2$ 
which converge to some $p\in I^2$, and let
$\rho:I^2\ra [1,1+c]$ be any continuous function
such that $\rho\restr_{S_k}=\rho_{k,\min(c,\frac{1}{k})}\circ\al_k$
where $\al_k:S_k\ra I^2$ is a similarity.

Also, note that to construct $\rho_{L,\bar c}$, we really
only need to construct a measurable function with 
the same property: if $\rho^k_{L,\bar c}$ is a sequence
of smoothings of a measurable function $\rho_{L,\bar c}$
which converge to $\rho_{L,\bar c}$ in $L^1$, then 
any sequence of $L$-biLipschitz maps
$\phi_k:I^2\ra \E^2$ with $Jac(\phi_k)=\rho^k_{L,\bar c}$
 will subconverge to a biLipschitz map $\phi:I^2\ra\E^2$
with $Jac(\phi)=\rho_{L,\bar c}$.

We now fix $L>1,\,c>0$, and explain how to construct 
$\rho_{L,c}$.  Let $R$ be the 
rectangle $[0,1]\times [0,\frac{1}{N}]\subset \E^2$,
where $N\gg 1$ is chosen suitably depending on $L$ and $c$,
and let $S_i=[\frac{i-1}{N},\frac{i}{N}]\times[0,\frac{1}{N}]$
be the $i^{th}$ square in $R$.
Define a ``checkerboard'' function $\rho_1:I^2\ra[1,1+c]$
by letting $\rho_1$ be $1+c$ on the squares $S_i$ with $i$ even 
and $1$ elsewhere.
Now subdivide $R$ into $M^2N$ squares using $M$ evenly spaced
horizontal lines and $MN$ evenly spaced vertical lines.
We call a pair of points {\em marked} if they are the
endpoints of a horizontal edge in the resulting grid.

The key step in the proof (Lemma \ref{MLE})
is to show that any biLipschitz
map $f:I^2\ra\E^2$ with Jacobian $\rho_1$ must stretch
apart a marked pair quantitatively more than it stretches
apart the pair $(0,0),\,(1,0)$; more precisely, there
is a $k>0$ (depending on $L,\, c$) so that 
$\frac{d(f(p),f(q))}{d(p,q)}>(1+k)d(f(0,0),f(1,0))$ for 
some marked pair $p,\,q$.  If this weren't true,
then we would have an $L$-biLipschitz map 
$f:I^2\ra\E^2$ which stretches apart all marked pairs by 
a factor of at most $(1+\eps)d(f(0,0),f(0,1))$, where $\eps\ll 1$.  This would
mean that
$f$ maps horizontal lines in $R$ to ``almost taut curves''.
Using triangle inequalities one checks that this
forces $f$ to map most marked pairs $p,\,q$
so that vector $f(q)-f(p)$
is $\approx d(p,q)(f(1,0)-f(0,0))$;
this in turn implies that for some $1\leq i\leq N$,
{\em all} marked pairs $p,\,q$ in the adjacent
squares $S_i,\,S_{i+1}$ 
are mapped by $f$ so that $f(q)-f(p)\approx d(p,q)
(f(1,0)-f(0,0))$.  Estimates then
show that $f(S_i)$ and $f(S_{i+1})$ have nearly the
same area, which 
contradicts the assumption that $Jac(f)=\rho_1$,
because $\rho_1$ is $1$ on one of the squares and $1+c$
on the other. 

Our next step is to  modify $\rho_1$ in a neighborhood 
of the grid in $R$:
we use thin rectangles (whose thickness will depend
on $L,\, c$) containing the horizontal edges of our grid,
and define $\rho_2:I^2\ra[1,1+c]$ by letting $\rho_2$
be a ``checkerboard'' function in each of these rectangles 
and $\rho_1$  elsewhere.
Arguing as in the previous paragraph, we will conclude
that some suitably  chosen pair of points in one of
these new rectangles will be stretched apart by a 
factor $>d(f(0,0),f(1,0))(1+k)^2$ under the map $f$.
Repeating this construction at smaller and smaller 
scales, we eventually obtain a function which 
can't be the Jacobian of an $L$-biLipschitz map.

The paper is organized as follows. In Section 2 we prove that Theorem
\ref{T2}
implies Theorem \ref{T1}. Section 3 is devoted to the proof of Theorem
\ref{T2}.

\section{Reduction of Theorem \ref{T1} to Theorem \ref{T2}}

Recall that every biLipschitz map is differentiable a.e., and the area
of
the image of a set is equal to the integral of the Jacobian over this
set.
We formulate our reduction as the following Lemma:

\begin{lemma}
 Let   $\rho:I^2\ra [1,1+c]$ be a measurable
function
which is not the Jacobian of any  biLipschitz
map $f:I^2\ra \E^2$ with 
\begin{equation}
Jac(f):=det(Df)=\rho \quad a.e.
\end{equation}
Then there is a separated net in $\E^2$ which is not
biLipschitz homeomorphic to $\Z^2$.
\end{lemma}

\proof
In what follows, the phrase ``subdivide the square
$S$ into subsquares will mean that $S$
is to be subdivided into squares using evenly 
spaced lines parallel to the sides of 
$S$.  Let ${\cal S}=\{S_k\}_{k=1}^\infty$ be a disjoint 
collection of square regions in $\E^2$ so that
each $S_k$ has integer vertices, sides parallel
to the coordinate axes, and the side length $l_k$
of $S_k$ tends to $\infty$ with $k$.  Choose
a sequence $m_k\in(1,\infty)$ with $\lim_{k\ra\infty}m_k=\infty$
and $\lim_{k\ra\infty}\frac{m_k}{l_k}=0$.
Let $\phi_k:I^2\ra S_k$ be the unique affine homeomorphism
with scalar linear part, and define $\rho_k:S_k\ra[1,1+c]$
by $\rho_k=(\frac{1}{\rho})\circ \phi_k^{-1}$. 
Subdivide $S_k$ into $m_k^2$ subsquares of side
length $\frac{l_k}{m_k}$.    Call this collection 
${\cal T}_k=\{T_{ki}\}_{i=1}^{m_k^2}$.
For each $i$ in $\{1,\ldots,m_k^2\}$, subdivide $T_{ki}$
into $n_{ki}^2$ subsquares $U_{kij}$ where $n_{ki}$
is the integer part of $\sqrt{\int_{T_{ki}}\rho_kd{\cal L}}$.
Now construct a separated net $X\subset\E^2$ by placing one
point at the center of each integer square
not contained in $\cup S_k$, and one point at the center
of each square $U_{kij}$.

We now prove the lemma by contradiction.  Suppose
$g:X\ra\Z^2$ is an $L$-biLipschitz homeomorphism.
Let $X_k=\phi_k^{-1}(X)\subset I^2$, and define
$f_k:X_k\ra\E^2$ by 
\begin{equation}
f_k(x)=\frac{1}{l_k}(g\circ\phi_k(x)-g\circ\phi_k(\star_k))
\end{equation}
where $\star_k$ is some basepoint in $X_k$.  Then
$f_k$ is an $L$-biLipschitz map from $X_k$ to 
a subset of $\E^2$, and the $f_k$'s are uniformly bounded.
By the proof of the Arzela-Ascoli theorem we may find a subsequence
of the $f_k$'s which ``converges uniformly'' to some biLipschitz
map $f:I^2\ra \E^2$.
By the construction of $X$, the counting measure on $X_k$
(normalized by the factor $\frac{1}{l_k^2}$)
converges weakly to $\frac{1}{\rho}$ times Lebesgue
measure, while the (normalized) counting measure on 
$f_k(X_k)$ converges weakly to Lebesgue measure.
It follows that $f_*((\frac{1}{\rho}){\cal L})={\cal L}\restr_{f(I^2)}$,
i.e. $Jac(f)=\rho$.
\qed

\section{Construction of a continuous function which is not a Jacobian
of a
biLipschitz map}

The purpose of this section is to prove Theorem \ref{T2}.
As explained in the introduction, Theorem \ref{T2} follows
from

\begin{lemma}
\label{ML}
For any given $L$ and $c>1$, there exists a continuous function
 $\rho:S=I^2\ra [1,1+c]$,
such that there is no $L$-biLipschitz 
homeomorphism $f:I^2\ra \E^2$ with 

$$Jac(f)=\rho \quad a.e.$$

\end{lemma}

\no{\em Proof of Lemma \ref{ML}.}
From now on, we fix two constants $L$ and $c$ and proceed to
construct of a  continuous 
function $\rho: I^2 \rightarrow [1,1+c]$ which is not a Jacobian of an
$L$-biLipschitz map.
By default, all functions which we will describe, take values between
$1$ and $1+c$. 

We say that two points $x,\,y\in I^2$ are $A$-stretched (under a map
$f:I^2 \rightarrow \E^2$)
if $d(f(x),f(y)) \geq Ad(x,y)$.

For $N\in\N$, $R_N$ be the rectangle $[0,1]\times [0,\frac{1}{N}]$ and
define a ``checkerboard'' function 
$\rho_N:R_N \rightarrow [1,1+c]$ by
$\rho_N(x,y)=1$ if $[Nx]$ is even and $1+c$ otherwise. It will be convenient
to introduce the squares 
$S_i=[\frac{i-1}{N},\frac{i}{N}] \times [0, \frac{1}{N}]$, $i=1, \dots,N$;  $\rho_N$ is
constant on the interior of each $S_i$.

The cornerstone of our construction is the following lemma:

\begin{lemma}
\label{MLE}
There are  $k>0,\,M,\,\mu$, and $N_0$ such that if 
$N\geq N_0$, $\eps\leq \frac{\mu}{N^2}$ 
then  the following holds: 
if the 
pair of points $(0,0)$
and $(1,0)$ is $A$-stretched under an $L$-biLipschitz
 map $f:R_N \rightarrow \E^2$
whose Jacobian differs from 
$\rho_N$ on a set of area no bigger than $\epsilon$ , then at least one
pair of points of the
form $((\frac{p}{NM}, \frac{s}{NM}), (\frac{q}{NM}, \frac{s}{NM}))$ 
is $(1+k)A$-stretched (where $p$
and 
$q$
are integers between 
$0$ and $NM$ and $s$ is an integer between  $0$ and $M$).

\end{lemma}

\no
{\em Proof of Lemma \ref{MLE}.}
We will need constants $k,\,l,\,m,\,\eps\in(0,\infty)$ 
and $M,\,N\in \N$, 
which will be chosen at the end of the argument. 
We will assume 
that $N>10$ and $c,\,l<1$.  Let $f:R_N\ra\E^2$ be an $L$-biLipschitz map
such that $Jac(f)=\rho_N$ off a set of measure $\eps$.  
Without loss of
generality we assume that $f(x)=(0,0)$ and $f(y)=(z,0)$, $z \geq A$.

We will use the notation $x_{pq}^i:=(\frac{p+M(i-1)}{NM}, \frac{q}{NM})$, where $i$ is 
an integer between 1 and $N$, and $p$ and $q$ are integers between 
$0$ and $M$.
We call these points {\em marked}. Note that the marked points in $S_i$
are precisely the vertices of the subdivision of
$S_i$ into $M^2$ subsquares.
The index $i$ gives the number of the square $S_i$, and $p$ and $q$ are
``coordinates'' of $x_{pq}^i$ within the square $S_i$.

We will prove Lemma \ref{MLE}  by
contradiction: we assume that all pairs of the form $x_{pq}^i, x_{sq}^j$ 
are no more than $(1+k)A$-stretched. 

If $x^i_{pq}\in S_i$ is a marked point, we say that $x^{i+1}_{pq}\in S_{i+1}$
is the {\em marked point corresponding to $x^i_{pq}$};  
corresponding points is obtained by adding the vector $(\frac{1}{N},0)$, 
where $\frac{1}{N}$ is the side length of the square $S_i$.
We are going to consider vectors between the images of marked points in
$S_i$ and the images of
corresponding marked points in the neighbor square $S_{i+1}$. 
We denote these vectors by $W_{pq}^i:=f(x_{pq}^{i+1})-f(x_{pq}^i)$.
We will see that most of the $W_{pq}^i$'s have to be extremely
close to the vector $W:=(A/N,0)$, and,  in particular, we will find a
square $S_i$ where $W_{pq}^i$ is
 extremely close to $W$ for {\em all} $0\leq p,q\leq M$. This will 
 mean that the areas of $f(S_i)$ and $f(S_{i+1})$ are very
close, since  $f(S_{i+1})$ is very close to a translate of
 $f(S_i)$. On the other hand,  except for a set of measure $\eps$,
 the 
Jacobian of $f$ is $1$ in one of the square $S_i,\,S_{i+1}$ and 
$1+c$ in the other. This allows us to estimate the 
difference of the
areas of their images from below and get a contradiction.

If $l\in (0,1)$, we say that a vector $W_{pq}^i=f(x_{pq}^{i+1})-f(x_{pq}^i)$, (or the 
marked point 
$x_{pq}^i$),  is {\em regular} if 
the length of its projection to the x-axis is
greater than $\frac{(1-l)A}{N}$. 
We say that a square $S_i$ is {\em regular} if all marked points  $x_{pq}^i$
in
this square are regular.

\no
{\em Claim 1.} There exist  $k_1=k_1(l)>0$, $N_1=N_1(l)$, 
such that if  $k \leq k_1,\,N\geq N_1$, there is a regular square.

\proof Reasoning by contradiction, we assume that all squares are
irregular. 
By the pigeon-hole principle, there is
a value of $ s $ (between 0 and $M$)  such that there are at least  
$\frac{N}{2M+2}$ irregular
vectors $W_{p_js}^{ i_j}$, $j=1,2,\dots J \geq \frac{N}{2M+2}$, where
 $i_j$ is an increasing sequence
with a fixed 
parity.  This means that we look for $l$-irregular vectors 
between marked points  in
the same row $s$ and only in squares $S_i$'s which have all indices
$i$'s even
or all odd.  The latter assumption guarantees that the segments
$[x_{p_js}^{ i_j},x_{p_js}^{ i_j+1}]$ do
not overlap.
We look at the polygon with  marked vertices
$$(0,0),x_{0s}^0= (0, s/MN),x_{p_1s}^{ i_1},
x_{p_1s}^{ i_1+1},x_{p_2s}^{ i_2},x_{p_2s}^{ i_2+1} \dots ,x_{p_Js}^{
i_J},
x_{p_Js}^{ i_J+1}, x_{Ms}^N=(1, s/MN), (1,0)$$
The  image of this polygon under $f$  connects $(0,0)$ and $(z,0)$ and, therefore, the
length of its projection
onto the x-axis is at least $z \geq A$. On the other hand, estimating
this
projection separately for
the images of $l$-irregular segments $[x_{p_js}^{ i_j},x_{p_js}^{ i_j+1}]$, 
the ``horizontal''
segments  $[x_{p_js}^{ i_j+1},x_{p_{j+1}s}^{ i_{j+1}}]$
and the  two ``vertical'' segments $[(0,0), x_{0s}^0]$ and $[x_{Ms}^N,
(1,0)]$  ,
one gets that the lengths of this projection is no bigger than

\begin{equation}
\label{projlength}
( \frac{N}{2M+2}) (\frac{(1-l)A}{N}) + 
(\frac{(1+k)A}{N})(N-\frac{N}{2M+2})+2\frac{L}{N}.
\end{equation}
\no
 The first term in (\ref{projlength}) bounds the total length of
projections of images of irregular segments
by the definition of irregular segments and total number of them. 
The second
summand is maximum stretch 
factor $(1+k)A$ between marked points times the total length of
remaining
horizontal segments. The third summand estimates the
lengths of images of segments   $[(0,0), x_{0s}^0]$ and $[x_{Ms}^N,
(1,0)]$
just by multiplying their lengths
by our fixed bound $L$ on the Lipschitz constant. 

Recalling that this projection is at least $z$, which in its turn is no
less than $A$, we get 
$$( \frac{N}{2M+2}) (\frac{(1-l)A}{N}) + 
(\frac{(1+k)A}{N})(N-\frac{N}{2M+2})+2\frac{L}{N}\geq A.$$
One easily checks that this is impossible when $k$ is sufficiently
small and $N$ is sufficiently large.  This 
contradiction
proves Claim 1.
\qed 

Let $W=(\frac{A}{N}, 0)$.

\no
{\em Claim 2.} Given any $m>0$, there is an $l_0=l_0(m)>0$ such that
if $l \leq l_0$ and $k \leq l$, then 
 $|W-W_{pq}^i| \leq \frac{m}{N}$ for every regular vector $W_{pq}^i$.

\proof Consider a regular vector $W_{pq}^i=(X,Y)$. 
 Since $W_{pq}^i$ is regular,  
$X \geq \frac{(1-l)A}{N}$. 
On the other hand, $X^2+Y^2 \leq \frac{(1+k)^2A^2}{N^2}$ and   $X \leq
\frac{(1+k)A}{N}$. Thus the difference of the 
$x$-coordinates of $W_{pq}^i$ and $W$ is bounded by $\frac{(l+k)A}{N} <\frac{2lA}{N}$.
Substituting the smallest possible value $\frac{(1-l)A}{N}$ for $X$
into 
 $X^2+Y^2 \leq \frac{(1+k)^2A^2}{N^2}$, we get 
 $Y^2 \leq \frac{2(l+k)A^2}{N^2} \leq
\frac{4lA^2}{N^2}$ . 
This implies that 

\begin{equation}
\label{w'sclose}
N|W-W_{pq}^i| \leq 2A\sqrt{l^2+l}\leq 2L\sqrt{l^2+l}.
\end{equation}   
The right-hand side of (\ref{w'sclose}) tends to zero with $l$,
so Claim 2 follows. 
\qed

\no
{\em Claim 3.}  There are $m_0>0$, $M_0$
such that if $m<m_0$ and $M>M_0$, then the 
following holds:
if for some $1\leq i\leq N$ and every $p,\,q$
we have $|W-W_{pq}^i|\leq\frac{m}{N}$, then
\begin{equation}
\label{areasclose}
|Area(f(S_{i+1}))-Area(f(S_i))|<\frac{c}{2N^2}.
\end{equation}

\proof
We assume that $i$ is even and therefore  $\rho$ takes the value $1$ on
$S_i$ and $1+c$ on $S_{i+1}$;
the other case is analogous.  We let $Q:=f(S_i)$ and 
$R=f(S_{i+1})$.

$Q$ is bounded by a curve (which is the image of the boundary of $S_i$).
Consider the result $\tilde R :=Q+W$ of translating $Q$ by the vector
$W=(A/N, 0)$.
The area of $\tilde R$ is equal to the area of $Q$.
 
The images of the marked points on the boundary of $S_i$ form
an $\frac{L}{NM}$-net on the 
boundary of $Q$, and the images
of marked points on the boundary of $S_{i+1}$ form an $\frac{L}{NM}$-net on 
$R$.
By assumption the difference between $W$ and each vector
$W_{pq}^i$
joining the image of a 
marked point on the boundary of $S_i$ and the 
image of the corresponding point on the
boundary  of $S_{i+1}$  is less than $\frac{m}{N}$.
We conclude that the boundary of $\tilde R$ lies within
the $\frac{m}{N}+\frac{2L}{MN}$-neighborhood of the boundary of
$R$.  Since $f$ is $L$-Lipschitz, the
length 
of the boundary of $R$ is $\leq 4L/N$. Using 
a standard estimate for the area of a neighborhood of a
curve, we obtain:
 
 $$|Area(R)-Area(Q)|=|Area(R)-Area(\tilde R)| 
\leq \frac{2L}{N}(\frac{m}{N}+\frac{2L}{MN})+\pi(\frac{m}{N}+\frac{2L}{MN})^2.$$
Therefore (\ref{areasclose}) holds if  $m$ is sufficiently small and 
 $M$ is  sufficiently large. 
 \qed
 
\no
{\em Proof of Lemma \ref{MLE} concluded.}
Now assume $m<m_0$, $M>M_0$, $l\leq l_0(m)$, $k\leq\min(l,k_1(l))$,
$N\geq N_1(l)$, and $\eps\leq\frac{c}{8N^2L^2}$.
Combining claims 1, 2, and 3, we find a square $S_i$
so that (\ref{areasclose}) holds.  On the other hand,
since $Jac(f)$ coincides with $\rho$ off a set of measure
$\eps$, $Area(f(S_i))\leq 1/N^2+\eps L^2$
and  $Area(f(S_{i+1})\geq (1+c)(1/N^2-\epsilon)$. Using
the assumption that  $\eps\leq\frac{c}{8N^2L^2}$
we get $$Area(R)-Area(Q)\geq \frac{c}{2N^2},$$
contradicting (\ref{areasclose}).
This contradiction proves Lemma \ref{MLE}.
\qed

\no
{\em Proof of  Lemma \ref{ML} continued.}
We will use an inductive construction based on Lemma \ref{MLE}.
Rather than dealing with an explicit construction of pairs of points
as in Lemma \ref{MLE}, it is more convenient
to us to use the following lemma, which is an obvious corollary of Lemma
\ref{MLE}. (To deduce this lemma from Lemma \ref{MLE}, just note
that all properties of interest persist if we scale our coordinate
system.)
 
\begin{lemma}
\label{IL}
There exists a constant $k>0$ such that, 
given any segment $\ol{xy}\subset I^2$ and any 
neighborhood $\ol{xy}\subset U\subset I^2$, 
there is a measurable function  $\rho:U \rightarrow [1,1+c]$, $\eps>0$
and a finite collection of 
non-intersecting segments $\ol{l_kr_k}\subset U$  with the
following property: 
if the pair $x,\,y$ is $A$-stretched by  an $L$-biLipschitz map $f:U
\rightarrow \E^2$   whose Jacobian differs from 
$\rho$ on a set of area  $<\epsilon$ , then 
for some $k$ the pair $l_k,\,r_k$ is $(1+k)A$-stretched
by $f$. 
The function  $\rho$ may be chosen to have finite image.
\end{lemma}

We will prove  Lemma \ref{ML} by induction, using 
the following statement.  (It is actually even slightly stronger than 
Lemma \ref{ML}  since it not only guarantees
non-existence of $L$-biLipschitz maps 
with a certain Jacobian, but also gives a finite collection of points,
such that at least one distance between 
them is distorted  more than by factor $L$.)

\begin{lemma}
\label{FF}
For each integer $i$  
there is a measurable function $\rho_i: I^2 \rightarrow [1, 1+c]$ , a finite
collection ${\cal S}_i$ of non-intersecting segments 
$\ol{l_kr_k}\subset I^2$, 
and $\epsilon_i>0$ with the following property: 
For every $L$-biLipschitz map $f:I^2 \rightarrow \E^2$   whose Jacobian
differs from 
$\rho_i$ on a set of area  $<\epsilon_i$ , at least one 
segment from ${\cal S}_i$  will have its endpoints
 $\frac{(1+k)^i}{L}$-stretched by $f$. 
\end{lemma}

\proof
The case $i=0$ is obvious.  Assume inductively that there are
$\rho_{i-1},\,\eps_{i-1}$, and a disjoint collection of
segments ${\cal S}_{i-1}=\{\ol{l_kr_k}\}$ which satisfy 
the conditions of the lemma.  Let $\{U_k\}$ be a disjoint
collection of open sets with $U_k\supset\ol{l_kr_k}$
and with total area $<\frac{\eps_{i-1}}{2}$.
For each $k$ apply Lemma \ref{IL} to $U_k$ to get a
function $\hat\rho_k:U_k\ra[1,1+c],\,\hat\eps_k>0$,
and a disjoint collection $\hat{\cal S}_k$ of segments.
Now let $\rho_i:I^2\ra[1,1+c]$ be the function which
equals $\hat\rho_k$ on each $U_k$ and equals $\rho_{i-1}$
on the complement of $\cup U_k$; let ${\cal S}_i=\cup \hat{\cal S}_k$,
and $\eps_i=\min\hat\eps_k$.  The required properties 
follow immediately.
\qed

Lemma \ref{ML} and (Theorem \ref{T2}) follows from Lemma \ref{FF}.

\bibliography{refs}
\bibliographystyle{alpha}

\end{document}